\documentclass[preprint,12pt]{elsarticle}

\usepackage{lipsum}
\usepackage{amsfonts}
\usepackage{graphicx}
\usepackage{epstopdf}
\usepackage{algorithmic}
\usepackage{amsmath}
\usepackage{multirow}
\usepackage{booktabs}
\usepackage{amsmath,amssymb}
\usepackage{units}
\usepackage{amssymb,amsmath}
\usepackage{ntheorem}
\usepackage{graphics,graphicx}
\usepackage{units}
\usepackage{color}
\usepackage{mathtools}
\usepackage{caption}

\newtheorem{definition}{Definition}[section]
\newtheorem{theorem}{Theorem}[section]

\newtheorem{lemma}{Lemma}[section]

\numberwithin{equation}{section}
\makeatletter

\newcommand{\Rmnum}[1]{\expandafter\@slowromancap\romannumeral #1@}
\makeatother
\journal{}

\begin{document}
	
	\begin{frontmatter}
		\title{A $\tau$ Matrix Based Approximate Inverse Preconditioning for Tempered Fractional Diffusion Equations
\footnote{This work was supported by the National Natural Science Foundation of China (No. 12001022).\\
$^\ast$ Corresponding author: Chaojie Wang.
\\
\indent {\emph{E-mail address: wangcj2019@btbu.edu.cn} }}
}
\author{Xuan Zhang$^{a}$, Chaojie Wang$^{a,\ast}$, Haiyu Liu$^{b}$}
\address{a.~School of Mathematics and Statistics, Beijing Technology and Business University, Beijing 100048, China}
\address{b.~School of Applied Science, Beijing Information Science and Technology University, Beijing 100192, China}
		\begin{abstract}
			Tempered fractional diffusion equations are a crucial class of equations widely applied in many physical fields. In this paper, the Crank-Nicolson method and the tempered weighted and shifts Grünwald formula are firstly applied to discretize the tempered fractional diffusion equations. We then obtain that the coefficient matrix of the discretized system has the structure of the sum of the identity matrix and a diagonal matrix multiplied by a symmetric positive definite(SPD) Toeplitz matrix. Based on the properties of SPD Toeplitz matrices, we use $\tau$ matrix approximate it and then propose a novel approximate inverse preconditioner to approximate the coefficient matrix. The $\tau$ matrix based approximate inverse preconditioner can be efficiently computed using the discrete sine transform(DST). In spectral analysis, the eigenvalues of the preconditioned coefficient matrix are clustered around 1, ensuring fast convergence of Krylov subspace methods with the new preconditioner. Finally, numerical experiments demonstrate the effectiveness of the proposed preconditioner.
		\end{abstract}
	\end{frontmatter}

	\section{Introduction}\label{section1}
	This paper discusses the numerical solution of variable coefficient tempered fractional diffusion equations (Tempered-FDEs) with initial-boundary value problems:
	\begin{equation}
		\left\{ \begin{array}{l}
			\frac{\partial u\left( x,t \right)}{\partial t}=d\left( x \right) \left( _a\boldsymbol{D}_{x}^{\beta ,\lambda}+_x\boldsymbol{D}_{b}^{\boldsymbol{\beta ,\lambda }} \right) u\left( x,t \right) +f\left( x,t \right) ,\\
			u\left( a,t \right) =0,\ u\left( b,t \right) =0,\ t\in \left[ 0,T \right] ,\,\,\\
			u\left( x,0 \right) =u_0\left( x \right) ,\ x\in \left[ a,b \right] ,\\
		\end{array} \right. 
	\end{equation}
	where \( f(x,t) \) is the source term, $d\left( x \right) \ge 0$ is the diffusion coefficient function, and $\lambda$ is a non-negative parameter. In Equation (1.1), $_a\boldsymbol{D}_{x}^{\beta ,\lambda}$ and $_x\boldsymbol{D}_{b}^{\beta ,\lambda}$ denotes the left and right Riemann-Liouville tempered fractional derivatives of the function \( u(x,t) \) with order $\beta \left( 1<\beta <2 \right)$, respectively defined by (see Baeumer and Meerschaert \cite{Baeumer2010}):
	\begin{equation}
		\nonumber
		_a\boldsymbol{D}_{x}^{\beta ,\lambda}u\left( x \right) =_aD_{x}^{\beta ,\lambda}u\left( x \right) -\beta \lambda ^{\beta -1}\partial _xu\left( x \right) -\lambda ^{\beta}u\left( x \right) 
	\end{equation}
	and
	\begin{equation}
		\nonumber
		_x\boldsymbol{D}_{b}^{\beta ,\lambda}u\left( x \right) =_xD_{b}^{\beta ,\lambda}u\left( x \right) +\beta \lambda ^{\beta -1}\partial _xu\left( x \right) -\lambda ^{\beta}u\left( x \right) 
	\end{equation}
	where
	\begin{equation}
		\nonumber
		_aD_{x}^{\beta ,\lambda}u\left( x \right) =\text{e}_{a}^{-\lambda x}D_{x}^{\beta}\left( \text{e}^{\lambda x}u\left( x \right) \right) =\frac{\text{e}^{-\lambda x}}{\Gamma \left( 2-\beta \right)}\frac{\partial ^2}{\partial x^2}\int_a^x{\frac{\text{e}^{\lambda s}u\left( s \right)}{\left( x-s \right) ^{\beta -1}}}ds
	\end{equation}
	and
	\begin{equation}
		\nonumber
		_xD_{b}^{\beta ,\lambda}u\left( x \right) =\text{e}_{x}^{\lambda x}D_{b}^{\beta}\left( \text{e}^{-\lambda x}u\left( x \right) \right) =\frac{\text{e}^{\lambda x}}{\Gamma \left( 2-\beta \right)}\frac{\partial ^2}{\partial x^2}\int_x^b{\frac{\text{e}^{-\lambda s}u\left( s \right)}{\left( s-x \right) ^{\beta -1}}}ds.
	\end{equation}
	
	Tempered fractional diffusion equations (Tempered-FDEs) are a crucial class of equations widely applied in fields such as biology, geophysics, and finance \cite{Magin2004}-\cite{Tang2023}. These equations are characterized by incorporating tempered fractional derivatives, modifying standard fractional diffusion equations by introducing an exponential tempering factor. This modification allows Tempered-FDEs to better simulate processes with finite propagation speed, accurately describing actual physical phenomena. However, analytical solutions to Tempered-FDEs and FDEs are generally difficult to obtain, leading to extensive research on the numerical solutions of FDEs in recent years \cite{Huang2021}-\cite{Shao2022}.
	
	For tempered fractional derivatives, stability of the resulting finite difference schemes can only be ensured if the spatial step size is sufficiently small when using standard fractional derivative approximations directly \cite{Wang2015}. A novel shifted approximation method for tempered fractional derivatives was proposed in \cite{Li2016}, successfully developing an unconditionally stable numerical scheme for one-dimensional (1D) tempered fractional diffusion equations with constant coefficients. Additionally, in \cite{Qu2017} proposed a Crank-Nicolson scheme for solving initial-boundary value problems of a class of variable coefficient tempered fractional diffusion equations, providing the discretization method used in this paper.
	
	Similar to fractional derivatives, tempered fractional derivatives also exhibit non-locality, resulting in dense or even full discrete coefficient matrices in the linear systems generated by discretization. Traditional direct solution methods, such as Gaussian elimination, require computational cost of \( O(N^3) \) and storage space of \( O(N^2) \), where \( N \) is the matrix size \cite{Wang2010}. Given the symmetric positive definite Toeplitz structure in the coefficient matrix, matrix-vector multiplications of the discrete coefficient matrix and the preconditioner can be efficiently implemented using fast Fourier transforms \cite{Chan1996, Huckle1994} and discrete sine transforms \cite{Bini1990,Lin2023}, respectively. For one-dimensional tempered space fractional diffusion equations, \cite{Tang2023} proposed a scaled diagonal and Toeplitz approximate splitting (SDTAS) preconditioning, theoretically proving that the spectra of the resulting preconditioned matrix cluster around 1. For one-dimensional space fractional diffusion equations, \cite{Zeng2022} introduced an approximate inverse preconditioning matrix based on the \(\tau\) matrix using the ideas of row-by-row inversion and interpolation, leveraging the properties of the symmetric positive definite Toeplitz matrix in the coefficient matrix. Results showed that preconditioned Krylov subspace methods converged at a superlinear rate with a total complexity of \( O(N \log N) \).
	
	In this paper, we use the second-order finite difference scheme proposed in \cite{Qu2017} to discretize the Tempered-FDEs. Specifically, we employ the Crank-Nicolson method for time discretization and the tempered weighted and shifts Grünwald formula for spatial discretization \cite{Li2016}. We keenly observe that the coefficient matrix of the resulting linear system is of the form \( I + DT \), where \( I \) is the identity matrix, \( D \) is a diagonal matrix, and \( T \) is a symmetric positive definite (SPD) Toeplitz matrix. Given the SPD property of the coefficient matrix, we approximate it using \(\tau\) matrix and construct a novel preconditioner using the row-by-row approximate inverse idea proposed by Pan and NG \cite{Pan2014}. Due to the SPD nature of the Toeplitz matrix in the coefficient matrix, the preconditioner can be efficiently computed using discrete sine transforms \cite{Bini1990}, with each step of the preconditioned Krylov subspace method having a computational complexity of \( O(N \log N) \). This construction method was first proposed in \cite{Zeng2022} for coefficient matrices of the form \( D + T \). Hence, we construct a row-by-row inverse preconditioner based on \(\tau\) matrix approximation and apply it for the first time to the coefficient matrix form \( I + DT \) generated by Tempered-FDEs. Theoretically, we prove that the spectra of the resulting preconditioned matrix cluster around 1, ensuring rapid convergence of Krylov subspace methods with the new preconditioner. Numerical experiments demonstrate the effectiveness of the proposed preconditioner, outperforming other preconditioners.
	
	The rest of this paper is organized as follows. In section 2, we present the discretized system from the Tempered-FDEs (1.1). In section 3, we construct an approximate inverse preconditioner. In section 4, We conduct spectral analysis of proposed preconditioned matrices. In section 5, numerical experiments are implemented to demonstrate the performance of the proposed preconditioner. Finally, we give a brief conclusion of the entire article in section 6.
	
	\section{Discretization of Tempered-FDEs}\label{section2}

	Let \( N \) and \( M \) be positive integers, representing the number of spatial and temporal partitions, respectively. We define the spatial step size \( h = \frac{a}{N+1} \) and the time step size \( \Delta t = \frac{T}{M} \). The spatial and temporal grids are given by \( x_i = ih \) for \( i = 0, 1, \ldots, N+1 \) and \( t_j = j\Delta t \) for \( j = 0, 1, \ldots, M \).
	
	Li and Deng \cite{Li2016} showed that the tempered fractional derivatives in (1.1) can be approximated using the tempered weighted and shifted Grünwald formula (Tempered-WSGD):
	\begin{equation}
		\nonumber
		\footnotesize
	_{a}D_{x}^{\beta ,\lambda}u\left( x_i,t_j \right) -\lambda ^{\beta}u\left( x_i,t_j \right) =\frac{1}{h^{\beta}}\sum_{k=0}^{i+1}{g}_{k}^{\left( \beta \right)}u\left( x_{i-k+1},t_j \right) -\frac{1}{h^{\beta}}\rho _{\beta}u\left( x_i,t_j \right) +O\left( h^2 \right) ,
	\end{equation}
	and
	\begin{equation}
		\footnotesize
		\nonumber
		_xD_{b}^{\beta ,\lambda}u\left( x_i,t_j \right) -\lambda ^{\beta}u\left( x_i,t_j \right) =\frac{1}{h^{\beta}}\sum_{k=0}^{N-i+2}{g}_{k}^{\left( \beta \right)}u\left( x_{i+k-1},t_j \right) -\frac{1}{h^{\beta}}\rho _{\beta}u\left( x_i,t_j \right) +O\left( h^2 \right) ,
	\end{equation}
	where $\rho _{\beta}=\left( \gamma _1e^{h\lambda}+\gamma _2+\gamma _3e^{-h\lambda} \right) \left( 1-e^{-h\lambda} \right) ^{\beta}$, the weights $g_{k}^{\left( \beta \right)}$ are given by
	\begin{equation}
		g_{k}^{\left( \beta \right)}=\begin{cases}
		\gamma _1w_{0}^{\left( \beta \right)}e^{h\lambda},&		k=0,\\
		\gamma _1w_{1}^{\left( \beta \right)}+\gamma _2w_{0}^{\left( \beta \right)},&		k=1,\\
		\left( \gamma _1w_{k}^{\left( \beta \right)}+\gamma _2w_{k-1}^{\left( \beta \right)}+\gamma 	_3w_{k-2}^{\left( \beta \right)} \right) e^{-\left( k-1 \right) h\lambda},&		k\geq 2.\\
		\end{cases}
	\end{equation}
	with $w_{0}^{\left( \beta \right)}=1,w_{k}^{\left( \beta \right)}=\left( 1-\frac{1+\beta}{k} \right) w_{k-1}^{\left( \beta \right)}, \text{ for }k\geq 1$ and $\gamma _1,\gamma _2,\gamma _3$ satisfy the linear system
	\begin{equation}
		\left\{ \begin{array}{l}
			\gamma _1+\gamma _2+\gamma _3=1,\\
			\gamma _1-\gamma _3=\beta /2.\\
		\end{array} \right. 
	\end{equation}
	
	Let $u_{i}^{j}\approx u\left( x_i,t_j \right) ,f_{i}^{j+\frac{1}{2}}=f\left( x_i,t_{j+\frac{1}{2}} \right) ,t_{j+\frac{1}{2}}=\frac{1}{2}\left( t_j+t_{j+1} \right) \ $ and $u_{i}^{j+\frac{1}{2}}=\frac{1}{2}\left( u_{i}^{j}+u_{i}^{j+1} \right) $ for $i=1,...,N,j=0,1,...,M\text{ and }1<\beta <2.$ At the mesh point $\left( x_i,t_{j+\frac{1}{2}} \right)$, We consider the Crank-Nicolson technique for time discretization and the tempered-WSGD approximation for the tempered fractional derivatives to discretize the tempered-FDEs in (1.1). Ignoring the truncation error, we obtain the following second-order finite difference scheme:
	\begin{equation}
		\footnotesize
		\begin{aligned}
			\frac{u_{i}^{j+1}-u_{i}^{j}}{\varDelta t}=&\frac{d_i}{2}\left( \frac{1}{h^{\beta}}\sum_{k=0}^{i+1}{g}_{k}^{\left( \beta \right)}u_{i-k+1}^{j+1}-\frac{1}{h^{\beta}}\rho _{\beta}u_{i}^{j+1}+\frac{1}{h^{\beta}}\sum_{k=0}^{N-i+2}{g}_{k}^{\left( \beta \right)}u_{i+k-1}^{j+1}-\frac{1}{h^{\beta}}\rho _{\beta}u_{i}^{j+1} \right)\\
			&+\frac{d_i}{2}\left( \frac{1}{h^{\beta}}\sum_{k=0}^{i+1}{g}_{k}^{\left( \beta \right)}u_{i-k+1}^{j}-\frac{1}{h^{\beta}}\rho _{\beta}u_{i}^{j}+\frac{1}{h^{\beta}}\sum_{k=0}^{N-i+2}{g}_{k}^{\left( \beta \right)}u_{i+k-1}^{j}-\frac{1}{h^{\beta}}\rho _{\beta}u_{i}^{j} \right) \\
			&+f_{i}^{j+\frac{1}{2}}.\\
		\end{aligned}
	\end{equation}

	We define the diagonal matrix as
	\begin{equation}
		\nonumber
		D=\text{diag}\left( d_1,d_2,...,d_n \right)
	\end{equation}
	and the vectors as
	\begin{equation}
		\nonumber
		u^j=\left[ u_{1}^{j},u_{2}^{j},...,u_{n}^{j} \right] ^T,f^{j+\frac{1}{2}}=\left[ f_{1}^{j+\frac{1}{2}},f_{2}^{j+\frac{1}{2}},...,f_{n}^{j+\frac{1}{2}} \right] ^T.
	\end{equation}
	
	For $j=0,1,...,m$ and a given $u^0$, The differential format (2.3) can be rewritten in the following matrix form:
	\begin{equation}
		\left( I+DG \right) u^{j+1}=\left( I-DG \right) u^j+\Delta tf^{j+\frac{1}{2}},
	\end{equation}
	where $I$ is the identity matrix, and $G$ is a symmetric positive definite Toeplitz matrix, defined as
	$G=G_{\beta}+G_{\beta}^{T}$, where $G_{\beta}^{T}$ represents the transpose of $G_{\beta}$. $G_{\beta}$ can be represented as
	\begin{equation}
		\nonumber
		G_{\beta}=-\frac{\varDelta t}{2h^{\beta}}\left( \begin{matrix}
			g_{1}^{\left( \beta \right)}&		g_{0}^{\left( \beta \right)}&		0&		...&		0&		0\\
			g_{2}^{\left( \beta \right)}&		g_{1}^{\left( \beta \right)}&		g_{0}^{\left( \beta \right)}&		0&		...&		0\\
			\vdots&		\ddots&		\ddots&		\ddots&		&		\vdots\\
			g_{N-1}^{\left( \beta \right)}&		\ddots&		\ddots&		\ddots&		&		g_{0}^{\left( \beta \right)}\\
			g_{N}^{\left( \beta \right)}&		g_{N-1}^{\left( \beta \right)}&		...&		...&		g_{2}^{\left( \beta \right)}&		g_{1}^{\left( \beta \right)}\\
		\end{matrix} \right)
	\end{equation} 
	
	 Hence, the coefficient matrix $A$ of Tempered-FDEs (2.3) is the sum of the identity matrix plus the diagonal matrix multiplied by a SPD Toeplitz matrix, which we denote as
	\begin{equation}
		A=I+DG
	\end{equation}

	By the following lemma 2.1, it is proven that the finite difference scheme (2.3) is unconditionally stable and second-order accurate in both space and time.
	
	\begin{lemma}
		
		Let S\ be the solution set of the linear system (2.2), if $1<\beta<2$
		$\lambda \geq 0,\left( \gamma _1,\gamma _2,\gamma _3 \right) \in \mathcal{S},and$
		\scriptsize
		\begin{flalign}
			\nonumber
			&\ \text{1.}\max \left\{ \frac{2\left( \beta ^2+3\beta -4 \right)}{\beta ^2+3\beta +2},\frac{\beta ^2+3\beta}{\beta ^2+3\beta +4} \right\} <\gamma _1<\frac{3\left( \beta ^2+3\beta -2 \right)}{2\left( \beta ^2+3\beta +2 \right)};or&\\
			\nonumber
			&\ \text{2.}\frac{\left( \beta -4 \right) \left( \beta ^2+3\beta +2 \right) +24}{2\left( \beta ^2+3\beta +2 \right)}<\gamma _2<\min \left\{ \frac{\left( \beta -2 \right) \left( \beta ^2+3\beta +4 \right) +16}{2\left( \beta ^2+3\beta +4 \right)},\frac{\left( \beta -6 \right) \left( \beta ^2+3\beta +2 \right) +48}{2\left( \beta ^2+3\beta +2 \right)} \right\} ;or&\\
			\nonumber
			&\ \text{3.}\max \left\{ \frac{\left( 2-\beta \right) \left( \beta ^2+\beta -8 \right)}{\beta ^2+3\beta +2},\frac{\left( 1-\beta \right) \left( \beta ^2+2\beta \right)}{2\left( \beta ^2+3\beta +4 \right)} \right\} <\gamma _3<\frac{\left( 2-\beta \right) \left( \beta ^2+2\beta -3 \right)}{2\left( \beta ^2+3\beta +2 \right)}&
		\end{flalign}
		\small
		we have
		$$
		g_{1}^{\left( \beta \right)}<0,\quad g_{2}^{\left( \beta \right)}+g_{0}^{\left( \beta \right)}>0,\quad g_{k}^{\left( \beta \right)}>0,\quad \forall k\geq 3,\quad \sum_{k=0}^{\infty}{g}_{k}^{\left( \beta \right)}=\rho _{\beta}\geq 0.
		$$
		
	\end{lemma}

	\section{Approximate Inverse Preconditioner Based on \(\tau\) Matrix}
	
	From the previous section, we know that at each time step, we need to solve a large linear system of equations (2.4), where the coefficient matrix is defined as in Equation (2.5). Since the coefficient matrix \(A\) is asymmetric, we can use Krylov subspace methods for solving, such as GMRES. To improve the performance and stability of Krylov subspace methods, preconditioning of the linear system is necessary. Pan and Ng\cite{Pan2014} proposed an approximate inverse preconditioner, and this preconditioning method has been widely used \cite{Baeumer2010, Chou2017, Du2019, Gan2023,Zhang2023}.
	
	Based on the fact $e_{i}^{\text{T}}A=e_{i}^{\text{T}}K_i$ and approximation $e_{i}^{\intercal}A^{-1}\approx e_{i}^{\intercal}K_{i}^{-1}$, where \({e}_i\) represents the \(i\)-th column of the identity matrix, $K_i=I+d_iG,i=1,2,\ldots,N.$ we obtain the approximate inverse preconditioner $P_1$:
	
	\begin{equation}
		\nonumber
		P_{1}^{-1}=\sum_{i=1}^N{e_ie_{i}^{T}K_{i}^{-1}},
	\end{equation}
	
	Since $K_i$ is a Toeplitz matrix, one method is to approximate it with a circulant matrix, yielding a circulant matrix based preconditioner:
	
	\begin{equation}
		\nonumber
		P_{C}^{-1}=\sum_{i=1}^N{e}_ie_{i}^{T}C_{i}^{-1},
	\end{equation}
	
	However, we note that $G$ is an SPD Toeplitz matrix, we can use the \(\tau\) matrix to approximate $G$\cite{Bini1990,Serra1999}. For convenience, the first column of the $G$ matrix can be expressed as $-\frac{\varDelta t}{2h^{\beta}}\left[ 2\left( g_{1}^{\left( \beta \right)}-\rho _{\alpha} \right) ,g_{0}^{\left( \beta \right)}+g_{2}^{\left( \beta \right)},g_{3}^{\left( \beta \right)},...,g_{N}^{\left( \beta \right)} \right] ^T$. According to \cite{Serra1999}, we can compute the \(\tau\) matrix approximation using Hankel correction:
	
	\begin{equation}
		\tau \left( G \right) =G-HC\left( G \right) ,
	\end{equation}
	where $HC\left( G \right)$ is a symmetric Hankel matrix, and its anti-diagonal elements can be computed from the corresponding elements in the following first column and last column:
	
	\begin{equation}
		-\frac{\varDelta t}{2h^{\beta}}\left( g_{3}^{\left( \beta \right)},g_{4}^{\left( \beta \right)},...,g_{N}^{\left( \beta \right)},0,0 \right) ^T,\quad -\frac{\varDelta t}{2h^{\beta}}\left( 0,0,g_{N}^{\left( \beta \right)},...,g_{4}^{\left( \beta \right)},g_{3}^{\left( \beta \right)} \right) ^T,	
	\end{equation}
	
	It is well-known that $\tau\left( G \right)$ defined in (3.1) can be diagonalized by the discrete sine transform\cite{Bini1990}:
	
	\begin{equation}
		\tau \left( G \right) =S\Lambda S,\quad S=\left( \sqrt{\frac{2}{N+1}}\cdot \sin \frac{\pi ij}{N+1} \right) ,\quad i,j=1,2,...,N,
	\end{equation}
	where $\Lambda =\text{diag}\left( \lambda _1,\lambda _2,...,\lambda _N \right)$, $S$ is the discrete sine transform matrix. 
	
	Clearly, $G_j=I+d_j\tau \left( G \right)$ remains a $\tau$ matrix. Hence, using the \(\tau\) matrix $G_j$ to approximate the SPD Toeplitz matrix $K_j$, we obtain the \(\tau\) matrix-based approximate preconditioner:
	
	\begin{equation}
		P_{2}^{-1}=\sum_{i=1}^N{e}_ie_{i}^{T}G_{i}^{-1}.
	\end{equation}
	
	Using the preconditioner $P_{2}^{-1}$, we need to compute $O\left( N \right)$ discrete Sine transforms for each iteration, which is still computationally intensive. Pan and Ng \cite{Pan2014} propose using interpolation methods can significantly reduce the computational complexity.
	
	Define $\{x_i\}_{i=1}^{N}$ as the set of all discrete points in the interval \([a, b]\), and the function $q_k\left( x \right) =\frac{1}{1+\lambda _kd\left( x \right)}$, where $\lambda _k\in sp\left( \tau \left( G \right) \right) =\{\lambda _1,\lambda _2,...,\lambda _N\},k=1,2,...,N$. Selects $l\left( l\ll N \right)$ interpolation points $\{\left( \tilde{x}_j,q_k\left( \tilde{x}_j \right) \right) \}_{j=1}^{l}$ in $\left( x_i,q_k\left( x_i \right) \right) \}_{i=1}^{N}$ and then uses piecewise linear interpolation to obtain the interpolation function: 
	
	\begin{equation}
		\begin{aligned}
			p_k\left( x \right) &=\phi _1\left( x \right) q_k\left( \tilde{x}_1 \right) +\phi _2\left( x \right) q_k\left( \tilde{x}_2 \right) +\cdots +\phi _l\left( x \right) q_k\left( \tilde{x}_l \right) \\
			&=\sum_{s=1}^l{\phi _s\left( x \right) q_k\left( \tilde{x}_s \right)}
		\end{aligned}			
	\end{equation}
	
	We use formula (3.3) to diagonalize $G_j$ as $G_j=S\Lambda _jS, j=1,2,...,N$, where $S$ is the discrete sine transform matrix, and $\varLambda _j$ is the diagonal matrix composed of the eigenvalues of $G_j$. Then, applying interpolation (3.5) to approximate $G_{j}^{-1}$, we have:
	
	\begin{equation}
		G_{j}^{-1}\approx S\left( \sum_{s=1}^l{\phi _s\left( x_j \right) \Lambda _{s}^{-1}} \right) S,
	\end{equation}
	
	Replacing $G_{j}^{-1}$ in (3.4) with approximation values (3.6), we obtain the following preconditioner:
	
	\begin{equation}
		\begin{aligned}
			P_{3}^{-1}&=\sum_{i=1}^N{e_ie_{i}^{T}S\left( \sum_{s=1}^l{\phi _s\left( x_i \right) \Lambda _{s}^{-1}} \right) S}\\
			&=\sum_{s=1}^l{\left( \sum_{i=1}^N{e_ie_{i}^{T}\phi _s\left( x_i \right)} \right)}S\left( \Lambda _{s}^{-1} \right) S\\
			&=\sum_{s=1}^l{\Phi _sS\Lambda _{s}^{-1}S},\\
		\end{aligned}
	\end{equation}
	where $\Phi _s=\text{diag}\left( \varphi _s\left( x_1 \right) ,\varphi _s\left( x_2 \right) ,...,\varphi _s\left( x_N \right) \right) $, $\Lambda _{s}^{-1}=\text{diag}\left( q_1\left( \tilde{x}_s \right) ,q_2\left( \tilde{x}_s \right) ,...,q_N\left( \tilde{x}_s \right) \right)$, $s=1,2,...,l$ are diagonal matrices. Finally, for a suitable number of interpolation points $l$, the computation of the preconditioner (3.7) requires \(O\left( lN\log N \right)\) operations, which is acceptable.
	
	\section{Spectral analysis}\label{section4}
	In this section, we discuss the spectral properties of the preconditioned matrix \( P_3^{-1}A \). We first present some off-diagonal decay property reasoning.
	 
	\begin{definition}(\cite{Stro2002})
		Let \(A = (a_{i,j})_{i,j \in I}\) be a matrix, where \(I = \mathbb{Z}, \mathbb{N}\) or \(\{1, 2, \ldots, N\}\), then we say \(A\) belongs to the class \(\mathcal{L}_s\), if 
		\[
		|a_{i,j}| \leq \frac{c}{(1 + |i - j|)^s}
		\]
		for \(s > 1\), and some constant \(c > 0\).
	\end{definition}
	
	\begin{lemma}(\cite{Meerschaert2008,Meerschaert2012})
		Let \( g_k^{(\beta)} \) be as defined in Section 2 with \( 1 < \beta < 2 \). Then it holds:
		\[
		g_k^{(\beta)} = (1 - \frac{\beta + 1}{k})g_{k-1}^{(\beta)}, \, k = 1, 2, \ldots,
		\]
		\[
		g_0^{(\beta)} = 1, \, g_1^{(\beta)} = -\beta < 0, \, 1 > g_2^{(\beta)} > g_3^{(\beta)} > \cdots > 0 \, \text{and}
		\]
		\[
		\sum_{k=0}^{\infty} g_k^{(\beta)} = 0, \, \sum_{k=0}^{m} g_k^{(\beta)} < 0, \, 1 \leq m < \infty.
		\]
		
	\end{lemma}
		
	\begin{lemma}(\cite{Wang2010})
		Suppose \( g_k^{(\beta)} \) is defined by Section 2 with \( 1 < \beta < 2 \). Then
		\[
		g_k^{(\beta)} = \frac{1}{\Gamma(-\beta)k^{\beta+1}} (1 + O(\frac{1}{k})),
		\]
		where \(\Gamma(x)\) is the Gamma function.
	\end{lemma}
	
	\begin{lemma}(\cite{Stro2002})
		Let \( G, A, K_i \) and \( T_\beta \) be as defined in Sections 2 and 3, then \( G^{-1}, A^{-1}, K_i^{-1} \) and \( G_\beta^{-1} \in \mathcal{L}_{\beta+1} \), i.e.,
		\[
		|L_{m,n}| \leq \frac{c_0}{(1 + |m - n|)^{\beta+1}},
		\]
		where \(c_0\) is a positive constant and \(L = (L_{m,n})_{m,n \in I}\) with \(I = \{1, 2, \ldots, N\}\) can be \(G, A, K_i,G_\beta, G^{-1}, A^{-1}, K_i^{-1} \) and \( G_\beta^{-1} \).
	\end{lemma}

	\begin{lemma}
		Let \(\beta \in (1, 2)\) and \(L \in \mathcal{L}_{\beta+1}\). Then \(\exists\) a constant \(\varpi\) s.t. \(\|L\|_{\infty} \leq \varpi\).
	\end{lemma}
	\textbf{Proof.} By making use of
	\[
	\sum_{k=q+1}^{\infty} \frac{1}{k^{\beta+1}} \leq \int_{q+1}^{\infty} \frac{1}{x^{\beta+1}} \, dx = \frac{1}{\beta q^{\beta}},
	\]
	we have 
	\[
	\sum_{k=q+1}^{\infty} \frac{1}{k^{\beta+1}} \leq \frac{1}{\beta}.
	\]
	Hence, 
	{\footnotesize
		\[
		\|L\|_{\infty} \leq \max_{1 \leq i \leq n} \sum_{j=1}^{n} |L_{i,j}| \leq \max_{1 \leq i \leq n} (|L_{i,i}| + \sum_{j \neq i} |L_{i,j}|) = \max_{1 \leq i \leq n} (c_0 + 2 \sum_{k=2}^{\infty} \frac{c_0}{k^{\beta+1}}) \leq \frac{(2+\beta)c_0}{\beta},
		\]}
	
	then letting \(\varpi = \frac{(2+\beta)c_0}{\beta}\), we obtain the result. \hfill \(\square\)

	Next, we discuss the spectral properties of \( P_3^{-1}A \) by analyzing the approximation \( P_3^{-1} - A^{-1} \). Since
	\[
	P_3^{-1} - A^{-1} = P_3^{-1} - P_2^{-1} + P_2^{-1} - P_1^{-1} + P_1^{-1} - A^{-1},
	\]
	we will separately discuss the approximations \( P_3^{-1} - P_2^{-1} \), \( P_2^{-1} - P_1^{-1} \), and \( P_1^{-1} - A^{-1} \).
	
	First, we analyze the property of \( P_1^{-1} - A \). We present Lemma 4.5 with detailed proofs available in \cite{Pan2014}
	
	\begin{lemma}
		$\textit{ for a given }\varepsilon > 0, \textit{ there exists a constant }c_{1}$ and $an$ integer $N_{1}$ such that for $l\geq N_{1}$ $we$ have
		
		$$\|P_1^{-1}-A^{-1}\|_\infty\leq c_1\max_{1\leq i\leq N}\Delta(x_i,l)+\varepsilon.$$
	where
		$$\Delta(x_i,l)=\max_{i-l<k<i+l}|x_k-x_i|=(l-1)h,$$
	\end{lemma}

	Furthermore, we discuss the property of $P_{2}^{-1}-P_{1}^{-1}$. By the definition of $P_1$ and $P_2$ in section 3, we have
	
	$$
	\begin{aligned}
		P_{2}^{-1}-P_{1}^{-1}&=\sum_{i=1}^N{e}_ie_{i}^{T}\left( G_{i}^{-1}-K_{i}^{-1} \right)\\
		&=\sum_{i=1}^N{e}_ie_{i}^{T}K_{i}^{-1}\left( K_i-G_i \right) G_{i}^{-1}\\
		&=\sum_{i=1}^N{d_ie}_ie_{i}^{T}K_{i}^{-1}\left( G-\tau \left( G \right) \right) G_{i}^{-1}.\\
	\end{aligned}
	$$
	
	We first prove through Lemma 4.6 that $G_i^{-1}$ is bounded.
	
	\begin{lemma}
		Let $G_{i}$ be defined in section 3. Then we have $\|G_{i}^{-1}\|_{\infty}<\eta^{-1}.$
	\end{lemma}
	\textbf{Proof}. The proof method can refer to Lemma 4.7 in \cite{Pan2014}, here $\eta=1$.
	
	As is shown in \cite{Bini1990} that $G-\tau(G)$ can be split as $E_1+F_1$ with $E_1$ being in accordance with $G-\tau(G)$ in the upper left and lower right $(N-1)\times(N-1)$ sub-matrices and vanishing in the other entries, and $rank(E_1)\leq2(N-1).$ Further, if we let $\varepsilon>0$ be fixed, then we have $\|F_1\|_2<\varepsilon.$

	Because $K_i^{-1}$ has the off-diagonal decay property, then $K_i^{-1}$ can be rewritten as $\tilde{K}_i+\hat{K}_i,$ where
	$$
	\tilde{K}_i=\left[ \begin{matrix}
		*&		\cdots&		*&		0&		\cdots&		0\\
		\vdots&		\ddots&		&		\ddots&		\ddots&		\vdots\\
		*&		&		\ddots&		&		\ddots&		0\\
		0&		\ddots&		&		\ddots&		&		*\\
		\vdots&		\ddots&		\ddots&		&		\ddots&		\vdots\\
		0&		\cdots&		0&		*&		\cdots&		*\\
	\end{matrix} \right] \quad \text{and\quad }\hat{K}_i=\left[ \begin{matrix}
		0&		\cdots&		0&		*&		\cdots&		*\\
		\vdots&		\ddots&		&		\ddots&		\ddots&		\vdots\\
		0&		&		\ddots&		&		\ddots&		*\\
		*&		\ddots&		&		\ddots&		&		0\\
		\vdots&		\ddots&		\ddots&		&		\ddots&		\vdots\\
		*&		\cdots&		*&		0&		\cdots&		0\\
	\end{matrix} \right] .
	$$
	
	here '*' denotes the nonzero entries.
		
	Considering
	\begin{align*}
		K_i^{-1}(G - \tau(G))G_i^{-1} &= K_i^{-1}(E_1 + F_1)G_i^{-1} \\
		&= K_i^{-1}E_1G_i^{-1} + K_i^{-1}F_1G_i^{-1} \\
		&= K_i^{-1}E_1G_i^{-1} + (\hat{K}_i + \hat{F}_i)F_1G_i^{-1} \\
		&= (K_i^{-1}E_1 + \hat{K}_iF_1)G_i^{-1} + \hat{F}_iF_1G_i^{-1},
	\end{align*}
	we derive
	\begin{align*}
		\|(K_{i}^{-1}E_{1}+\hat{K}_{i}F_{1})G_{i}^{-1}\|_{\infty}& \leq(\|K_i^{-1}\|_\infty\cdot\|E_1\|_\infty+\|\hat{K}_i\|_\infty\cdot\|F_1\|_\infty)\cdot\|G_i^{-1}\|_\infty  \\
		&\leq\varepsilon(\|K_{i}^{-1}\|_{\infty}+\|F_{1}\|_{\infty})\cdot\|G_{i}^{-1}\|_{\infty}
	\end{align*}
	and
	\begin{align*}
		\parallel\sum_{i=1}^{N}e_{i}e_{i}^{G}(K_{i}^{-1}E_{1}+\hat{K}_{i}F_{1})G_{i}^{-1}\parallel_{\infty}& =\max_{1\leq i\leq N}\|e_{i}^{G}(K_{i}^{-1}E_{1}+\hat{K}_{i}F_{1})G_{i}^{-1}\|_{1} \\
		&=\max_{1\leq i\leq N}\|(K_{i}^{-1}E_{1}+\hat{K}_{i}F_{1})G_{i}^{-1}\|_{\infty} \\
		&\leq\varepsilon\cdot\max_{1\leq i\leq N}(\|K_{i}^{-1}\|_{\infty}+\|F_{1}\|_{\infty})\cdot\|G_{i}^{-1}\|_{\infty}.
	\end{align*}
	Since $\|K_i^{-1}\|_{\infty}$ and $\|F_1\|_{\infty}$ are limited due to the off-diagonal decay characteristics of $K_i^{-1}$ and $G$, it follows that $\|G_i^{-1}\|_{\infty}$ is also bounded. In other words, there is a constant $c_2 > 0$ such that
	\[
	\left\|\sum_{i=1}^{N} e_i^T (K_i^{-1}E_1 + \hat{K}_iF_1)G_i^{-1} \right\|_1 \leq c_2 \cdot \epsilon.
	\]
	
	Next, we focus on the matrix product $\hat{K}_iF_1G_i^{-1}$. Suppose the dimension of the blocks $\hat{K}_i$ matches the dimension of the block $F_1$. Otherwise, the smaller one can be extended. By examining the structure of $\hat{K}_i$ and $F_1$, we can compute the following results:
	$$\tilde{K}_{i}F_{1}G_{i}^{-1}=\begin{pmatrix}0&0&\ldots&\ldots&0&+\\0&0&\ldots&\ldots&0&+\\0&0&\ldots&\ldots&0&0\\\vdots&\vdots&&&\vdots&\vdots\\0&0&\ldots&\ldots&0&0\\+&0&\ldots&\ldots&0&0\\+&0&\ldots&\ldots&0&0\end{pmatrix}, G_{i}^{-1}=\begin{pmatrix}+&+&\ldots&\ldots&+&+\\+&+&\ldots&\ldots&+&+\\0&0&\ldots&\ldots&0&0\\\vdots&\vdots&&&\vdots&\vdots\\0&0&\ldots&\ldots&0&0\\+&+&\ldots&\ldots&+&+\\+&+&\ldots&\ldots&+&+\end{pmatrix},$$
	where ‘+’ denotes nonnegative entries. Therefore, we have
	\[
	\text{rank}\left( \sum_{i=1}^{N} e_i^T \hat{K}_iF_1G_i^{-1} \right) \leq 4\varsigma,
	\]
	where $\varsigma$ is the dimension of the blocks in $\hat{K}_iF_1G_i^{-1}$. Consequently,
	\[
	\sum_{i=1}^{N} e_i^T \left( K_i^{-1} - K_i^{-1} \right) = \sum_{i=1}^{N} e_i^T \left( K_i^{-1}E_1 + \hat{K}_iF_1 \right)G_i^{-1} + \sum_{i=1}^{N} e_i^T \hat{F}_iF_1G_i^{-1},
	\]
	is a sum of a small norm matrix and a low rank matrix. To summarize, we obtain the following results.
	
	\begin{theorem}
		Let $P_1$ and $P_2$ be defined in Section 4. The approximation $P_2^{-1}$ to $P_1^{-1}$ satisfies
		\[
		P_2^{-1} - P_1^{-1} = E_2 + F_2,
		\]
		where $E_2$ and $F_2$ are matrices of small norm and low rank, respectively, i.e., $\|E_2\| < c_2 \cdot \epsilon$ and $\text{rank}(F_2) \leq 4\varsigma$.
	\end{theorem}

	We now consider the approximation of the preconditioner $P_2$ and the coefficient matrix $A$.
	
	\begin{theorem}
		Let $P_1$, $P_2$, and $A$ be defined as previously stated. Then there exists $N_2$ such that for $N > N_2$, we have
		$$P_{2}^{-1}-A^{-1}=E_{P_{2}}+F_{P_{2}},$$
		where $E_{P_{2}}$ and $ F_{P_{2}}$ are of small norm and of low rank, respectively.
	\end{theorem}
	
	\textbf{Proof.} Due to Lemma 4.5, we know
	$$\|P_2^{-1}-A^{-1}\|_\infty\leq\max_{1\leq t\leq N}\Delta(x_i,N_1)+\varepsilon=c_1(N_1-1)h+\varepsilon.$$
	$\text{Let } N_2 \text{ be an integer such that } (N_1-1)h<\varepsilon\text{, then}$
	$$\begin{aligned}
		P_{2}^{-1}-A^{-1}& =P_{2}^{-1}-P_{1}^{-1}+P_{1}^{-1}-A^{-1} \\
		&=E_2+F_2+P_1^{-1}-A^{-1} \\
		&=E_2+(P_1^{-1}-A^{-1})+F_2 \\
		&\triangleq E_{P_{2}}+F_{P_{2}}, \\
	\end{aligned}$$
	where $E_{P_{2}}=E_{2}+(P_{1}^{-1}-A^{-1})$ and $F_{P_{2}}=F_{2}.$ As
	$$\|E_{P_2}\|_\infty=\|E_2\|_\infty+\|P_1^{-1}-A^{-1}\|_\infty<(c_1+c_2+1)\cdot\varepsilon,$$
	then the conclusion follows.$~\square $
	
	Finally, we discuss the properties of $P_{3}^{-1}-P_{2}^{-1}$. We present the following conclusion, with detailed proof provided in \cite{Zeng2022,Pan2014}.
	\begin{theorem}
		Suppose $\varpi$ is sufficiently small and $l \ll N$. Denote by $\varpi = \max_{1 \leq i \leq N} \max_{1 \leq j \leq N} \{ | p_{\lambda_j}(d_i) - q_{\lambda_j}(d_i) | \}$. Then for any given $\epsilon > 0$, $\exists N_3 > 0$ (independent of $N$) such that
		\[
		P_3^{-1} - P_2^{-1} = E_3 + F_3,
		\]
		where $E_3$ and $F_3$ satisfy $\|E_3\|_{\infty} \leq \varpi (2N_3 + 1) + \epsilon$ and $F_3$ is of a low rank matrix, respectively.
	\end{theorem}
	By combining the above theorems, we arrive at the following conclusion.
	
	\begin{theorem}
		$Let P_{3} and$ A be defined previously, then $\exists$ an integer N$_3,such$ that for N $>N_3$ it holds
		$$P_3^{-1}-A^{-1}=E_{P_3}+F_{P_3},$$
		where E$_{P_{3}}$ and F$_{P_{3}}$ are of small norm and of low rank, respectively.
	\end{theorem}

	In the following section, we provide numerical experiments to illustrate the effectiveness of the proposed preconditioner.

	\section{Numerical experiments}\label{section5}
	In this section, we conduct numerical experiments on the studied Tempered-FDEs (1.1) to verify the efficiency of the proposed preconditioning matrices $P_{3}^{-1}$ (3.7) denoted as $P_{TAI}^{-1}\left( l \right)$ , where $l$ represents the number of interpolation points. For comparison, we also use the Strang circulant approximate inverse matrix, the diagonal Toeplitz splitting preconditioning matrix, and the skew-diagonal Toeplitz splitting preconditioning matrix based on $\tau$ matrices. These three preconditioning matrices are denoted as $P_{CAI}^{-1}\left( l \right)$, $P_{DCS}^{-1}$ and $P_{SDTAS_{\tau}}^{-1}$, respectively. All experiments are  performed via Matlab (version R2022b) on a personal computer with a 4.00GHz AMD Ryzen 9 7940H CPU, 16.00 GB of memory, and the Windows 11 operating system.
	
	In all experiments, we searched for the optimal number of interpolation points $l$ in the proposed preconditioning matrices $P_{TAI}^{-1}\left( l \right)$. Therefore, for two experiments, we choose the optimal number of interpolation points $l$ to be 8 and 12, respectively. We specify the same grid density in both space and time $N=\varDelta t=h$, and for a certain initial point $u^{\left( 0 \right)}$, we define the relative error of the $k$-th iteration as follows:
	$$
	RES=\frac{\lVert b-Au^{\left( k \right)} \rVert _2}{\lVert b-Au^{\left( 0 \right)} \rVert _2}
	$$
	Let the initial point be a zero vector, the relative error tolerance be $RES<10^{-7}$, and the maximum number of iterations be 1000. We solve the problem using GMRES with the above preconditioners.
	
	For the following variable-coefficient tempered fractional diffusion equations (Tempered-FDEs):
	\begin{equation}
		\begin{cases}\frac{\partial u(x,t)}{\partial t}=d(x)\left(_0\mathbf{D}_x^{\beta,\lambda}+_x\mathbf{D}_1^{\beta,\lambda}\right)u(x,t)+f(x,t),\\u(0,t)=0, u(1,t)=0, t\in[0,1],\\u(x,0)=0, x\in[0,1].\end{cases}
	\end{equation}
	The exact solution of the equation (5.1) is $u\left( x,t \right) \simeq te^{-\lambda x}x^3\left( 1-x \right) ^3$. The source term is given as follows:
	$$
	\begin{array}{c}
		f\left( x,t \right) = e^{-\lambda x}x^3\left( 1-x \right) ^3-td\left( x \right) \left[ e^{-\lambda x}\sum_{m = 0}^3{\left( \left( -1 \right) ^m\left( \begin{array}{c}
				3\\
				m\\
			\end{array} \right) \frac{\Gamma \left( 4+m \right)}{\Gamma \left( 4+m-\beta \right)}x^{3+m-\beta} \right)} \right.\\
		\left. +e^{\lambda \left( x-2 \right)}\sum_{j = 0}^{30}{\frac{3^j}{j!}}\left( \sum_{m = 0}^3{\left( \left( -1 \right) ^m\left( \begin{array}{c}
				3\\
				m\\
			\end{array} \right) \frac{\Gamma \left( 4+m+j \right)}{\Gamma \left( 4+m+j-\beta \right)}\left( 1-x \right) ^{j+3+m-\beta} \right)} \right) \right]\\
		+2td\left( x \right) \left( \lambda ^{\beta} \right) e^{-\lambda x}x^3\left( 1-x \right) ^3\\
	\end{array}
	$$
	
	In this example, we conduct numerical experiments using the following two coefficients $d\left( x \right)$ of different complexity:
	$$
	d_1\left( x \right) =\frac{e^{5x}}{1+x}
	$$
	
	and 
	$$
	d_2\left( x \right) =\frac{e^{3x}+0.2}{x\left( 1-x \right)}
	$$
	
	The coefficient $d_1\left( x \right)$ and $d_2\left( x \right)$ is continuous over $\left[ 0,1 \right]$ but $d_2\left( x \right)$ has a singularity at point $x=0$ and $x=1$. Below, we present the numerical experiments for different coefficients of the equation.
	
	In the first experiment, We choose the coefficients $d_1\left( x \right) =\frac{e^{5x}}{1+x}$, $\lambda =1.5$, $\beta =1.2$ and the remaining parameters $\gamma _2$, $\gamma _3$ are computed from the linear system (2.2). To illustrate the effectiveness of our proposed preconditioning matrices, we compared $P_{TAI}^{-1}\left( l \right)$, $P_{CAI}^{-1}\left( l \right)$, $P_{DCS}^{-1}$ and $P_{SDTAS_{\tau}}^{-1}$. The parameters for the first two preconditioning matrices are $\l=8$, and the parameters for the latter two are $\beta =0.10$.

	we plot the two-dimensional images of the analytical solution and the numerical solution obtained using the preconditioners, as shown in Fig. 5.1:
	
	\begin{figure}[!htbp]
		\setlength{\abovecaptionskip}{-0cm}
		\setlength{\belowcaptionskip}{-0.25cm}
		\centering
		\includegraphics[height=5cm]{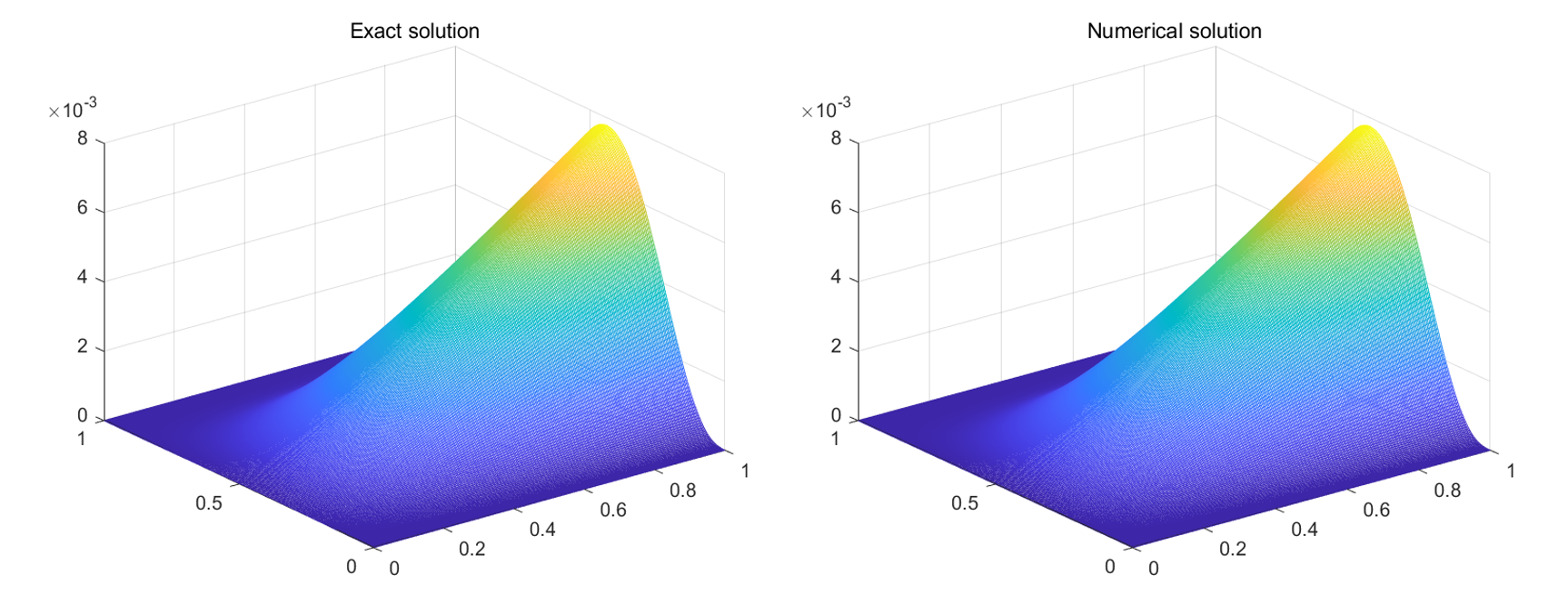}
		\caption*{Fig. 5.1. Comparison of the exact solutions and $P_{TAI}^{-1}\left( 10 \right)$ numerical solutions($d=d_1\left( x \right)$, $\lambda =1.5$, $\beta =1.2$, $\gamma _1=0.75$)}
		\label{figure1}
	\end{figure}
	
	We use GMRES with $P_{TAI}^{-1}\left( l \right)$, $P_{CAI}^{-1}\left( l \right)$, $P_{DCS}^{-1}$ and $P_{SDTAS_{\tau}}^{-1}$ preconditioners as well as GMRES without preconditioning to test. The results of the computations at discretization levels $N=10^8,10^9,10^{10},10^{11},10^{12}$ are shown in Table 1. Here, "IT" denotes the average number of iterations per step, and "CPU" denotes the total computation time (s) of the algorithm. To facilitate the comparison of the performance of different preconditioning matrices, we visualized the results from Table 1, as shown in Fig. 5.2. From Table 1 and Fig. 5.2, we can see that the performance of $P_{TAI}^{-1}\left( l \right)$ is significantly better than that of $P_{CAI}^{-1}\left( l \right)$, $P_{DCS}^{-1}$ and $P_{SDTAS_{\tau}}^{-1}$. The number of iterations for each preconditioning method remains relatively stable as the discretization level increases, but the number of iterations for $P_{TAI}^{-1}\left( l \right)$ is notably lower than that for the other three. Similarly, in terms of computation time, the computation time for $P_{TAI}^{-1}\left( l \right)$ is much lower than that for the other three, and the growth in computation time for $P_{TAI}^{-1}\left( l \right)$ is slower as the discretization level increases.
	
	\begin{table}[htbp]
		\centering
		\caption{Numerical results of $P_{TAI}^{-1}\left( l \right)$, $P_{CAI}^{-1}\left( l \right)$, $P_{DCS}^{-1}$, $P_{SDTAS_{\tau}}^{-1}$ and $I$($d=d_1\left( x \right)$, $\lambda =1.5$, $\beta =1.2$, $\gamma _1=0.75$)}
		\begin{tabular}{ccccccc}
			\toprule
			\multirow{2}[4]{*}{\textbf{Method}} & \multirow{2}[4]{*}{\textbf{Index}} & \multicolumn{5}{c}{\textbf{N}} \\
			\cmidrule{3-7}          &       & $2^8$   & $2^9 $    & $2^{10}$    & $2^{11}$     &$ 2^{12} $\\
			\midrule
			\multirow{3}[1]{*}{$P_{TAI}^{-1}\left( l \right)$} & $\l$     & 8     & 8     & 8     & 8     & 8 \\
			& IT    & \textbf{6.02 } & \textbf{5.01 } & \textbf{5.00 } & \textbf{5.00 } & \textbf{4.00 } \\
			& CPU   & \textbf{0.19 } & \textbf{0.46 } & \textbf{1.49 } & \textbf{5.87 } & \textbf{24.03 } \\
			\multirow{3}[0]{*}{$P_{CAI}^{-1}\left( l \right)$} & $\l$     & 8     & 8     & 8     & 8     & 8 \\
			& IT    & 8.03  & 7.01  & 6.01  & 6.00  & 6.00  \\
			& CPU   & 0.25  & 0.79  & 2.49  & 9.99  & 37.84  \\
			\multirow{3}[0]{*}{$P_{DCS}^{-1}$} & $\beta$     & 0.1   & 0.1   & 0.1   & 0.1   & 0.1 \\
			& IT    & 16.06  & 16.03  & 16.02  & 16.01  & 17.00  \\
			& CPU   & 0.38  & 1.31  & 4.53  & 17.33  & 65.74  \\
			\multirow{3}[0]{*}{$P_{SDTAS_{\tau}}^{-1}$} & $\beta$     & 0.1   & 0.1   & 0.1   & 0.1   & 0.1 \\
			& IT    & 11.04  & 12.02  & 14.01  & 15.01  & 16.00  \\
			& CPU   & 0.20  & 0.52  & 1.91  & 7.32  & 36.86  \\
			\multirow{2}[1]{*}{$I$} & IT    & 87.34 & 108.21 & 123.12 & 134.07 & 140.03 \\
			& CPU   & 0.82  & 3.66  & 14.33 & 70.64 & 555.72 \\
			\bottomrule
		\end{tabular}%
		\label{tab:addlabel}%
	\end{table}%

	\begin{figure}[!htbp]
		\setlength{\abovecaptionskip}{-0cm}
		\setlength{\belowcaptionskip}{-0.25cm}
		\centering
		\includegraphics[height=5.5cm]{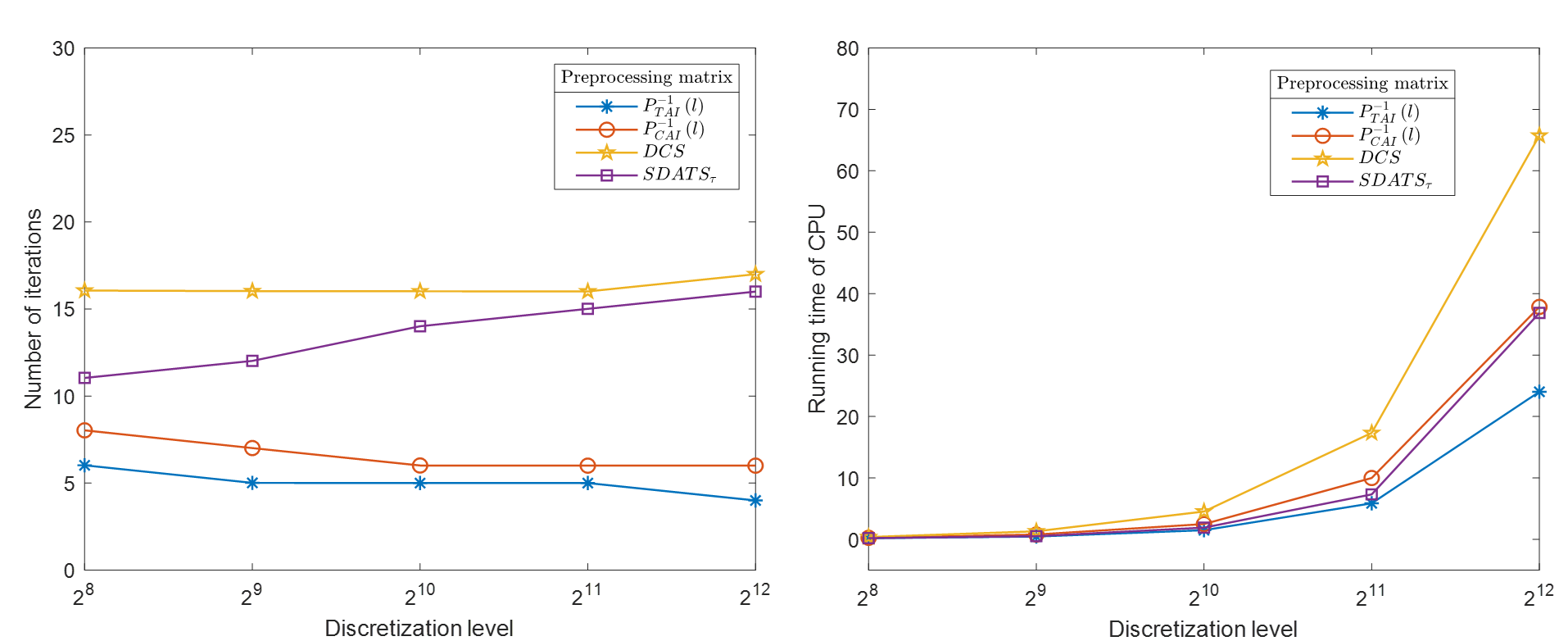}
		\caption*{Fig. 5.2.  Comparison of average iterations and computation time for $P_{TAI}^{-1}\left( l \right)$, $P_{CAI}^{-1}\left( l \right)$, $P_{DCS}^{-1}$, $P_{SDTAS_{\tau}}^{-1}$($d=d_1\left( x \right)$, $\lambda =1.5$, $\beta =1.2$, $\gamma _1=0.75$)}
		\label{figure2}
	\end{figure}
	
	 To further illustrate the effectiveness of the proposed preconditioning matrices, we plot the eigenvalue distributions of the preconditioned coefficient matrices $P_{TAI}^{-1}\left( l \right)$, $P_{CAI}^{-1}\left( l \right)$, $P_{DCS}^{-1}$, $P_{SDTAS_{\tau}}^{-1}$ and the original coefficient matrix $A$ when $N=10^8$, as shown in Fig. 5.3. It can be observed that the eigenvalue distribution of the coefficient matrix $A$ without preconditioning is highly dispersed, while the eigenvalue distributions of the preconditioned matrices $AP_{TAI}^{-1}\left( l \right)$ are concentrated around 1. Notably, the eigenvalue distribution of $AP_{TAI}^{-1}\left( l \right)$ is more concentrated than those of the other preconditioners.
	
	\begin{figure}[!htbp]
		\setlength{\abovecaptionskip}{-0cm}
		\setlength{\belowcaptionskip}{-0.25cm}
		\centering
		\includegraphics[height=7.5cm]{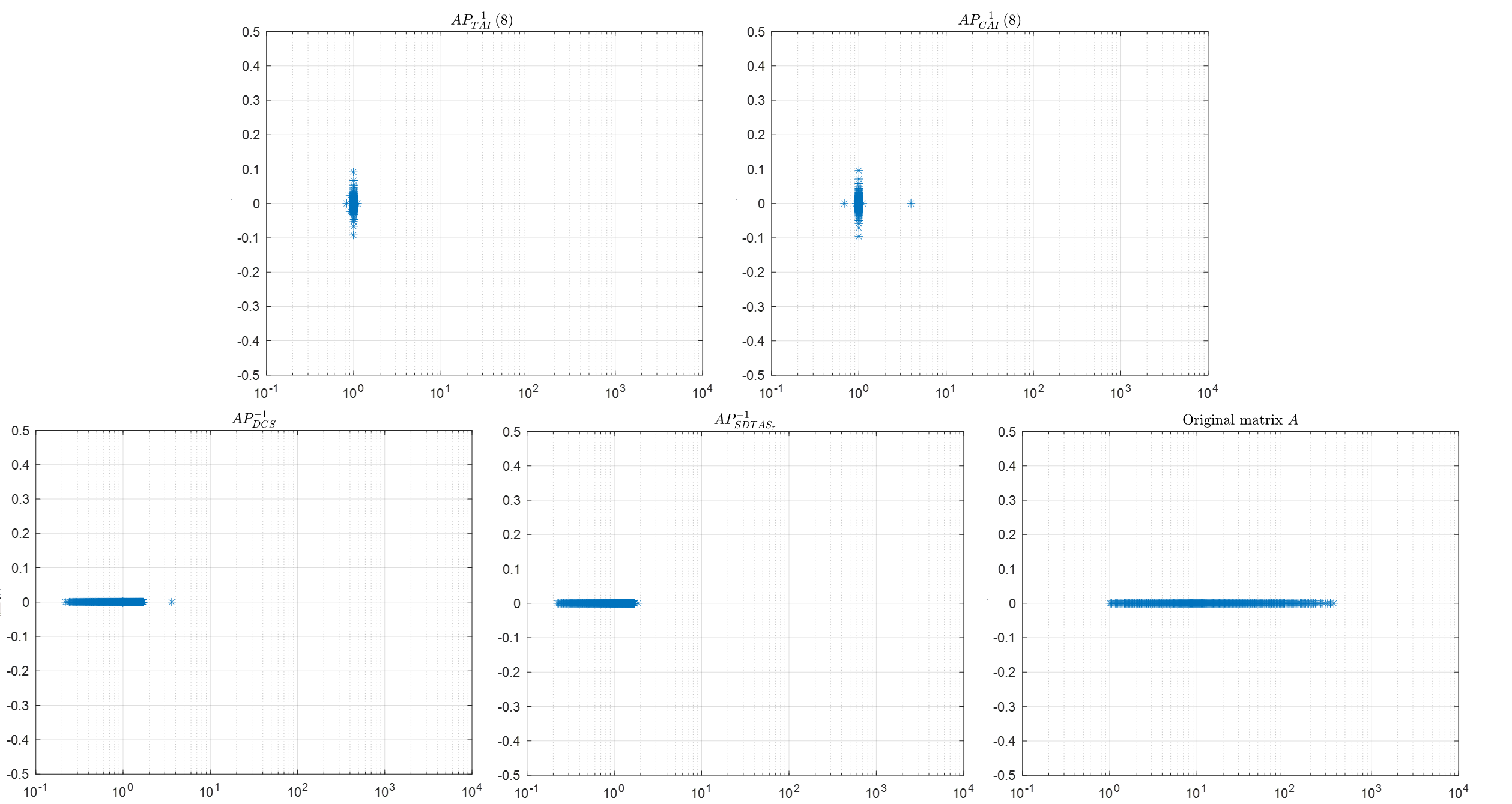}
		\caption*{Fig. 5.3.  The eigenvalue distributions of the preconditioned matrices $AP_{TAI}^{-1}\left( l \right)$, $AP_{CAI}^{-1}\left( l \right)$, $AP_{DCS}^{-1}$, $AP_{SDTAS_{\tau}}^{-1}$ and the original matrix $A$($d=d_1\left( x \right)$, $\lambda =1.5$, $\beta =1.2$, $\gamma _1=0.75$)}
		\label{figure3}
	\end{figure}
	
	In the second experiment, we chose the coefficient $d=d_2\left( x \right)$, which has a singularity at point $x=0$ and $x=1$, and keep other parameters unchanged. Similarly, to illustrate the effectiveness of our proposed preconditioning matrices, we compared $P_{TAI}^{-1}\left( l \right)$, $P_{CAI}^{-1}\left( l \right)$, $P_{DCS}^{-1}$ and $P_{SDTAS_{\tau}}^{-1}$. The parameters for the first two preconditioning matrices were $l=12$ and the parameters for the latter two were $\beta =0.15$.
	
	We plot the two-dimensional images of the excat solution and the numerical solution $P_{TAI}^{-1}\left( l \right)$ obtained using the preconditioners, as shown in Fig. 5.4. We use GMRES with $P_{TAI}^{-1}\left( l \right)$, $P_{CAI}^{-1}\left( l \right)$, $P_{DCS}^{-1}$ and $P_{SDTAS_{\tau}}^{-1}$ preconditioners as well as GMRES without preconditioning to test. The results of the computations at discretization levels $N=10^8,10^9,10^{10},10^{11},10^{12}$ are shown in Table 2 and visualized in Fig. 5.5. In terms of the number of iterations, all preconditioning methods show relatively stable iteration numbers as the discretization level increases, but the number of iterations for $P_{TAI}^{-1}\left( l \right)$ is notably lower than that for the other three. In terms of computation time, the computation time for $P_{TAI}^{-1}\left( l \right)$ is slightly higher than that for $P_{SDTAS_{\tau}}^{-1}$ at discretization level  $N=10^8,10^9,10^{10}$, but the difference is not significant. As the discretization level increases, the computation time for $P_{TAI}^{-1}\left( l \right)$ grows more slowly and is significantly lower than that for the other three preconditioners.
	
	\begin{figure}[!htbp]
		\setlength{\abovecaptionskip}{-0cm}
		\setlength{\belowcaptionskip}{-0.25cm}
		\centering
		\includegraphics[height=5cm]{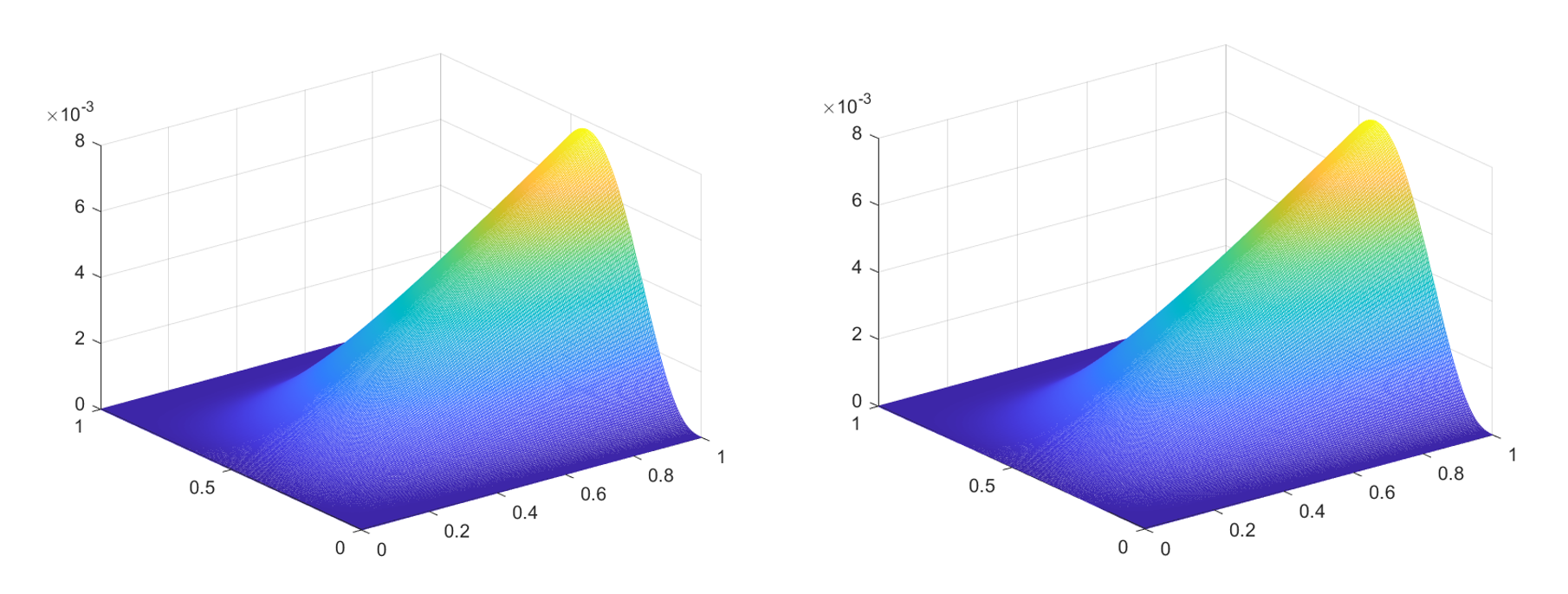}
		\caption*{Fig. 5.4. Comparison of the exact solutions and $P_{TAI}^{-1}\left( l \right)$ numerical solutions($d=d_2\left( x \right)$, $\lambda =1.5$, $\beta =1.2$, $\gamma _1=0.75$)}
		\label{figure4}
	\end{figure}
	
	\begin{table}[htbp]
		\centering
		\caption{Numerical results of $P_{TAI}^{-1}\left( l \right)$, $P_{CAI}^{-1}\left( l \right)$, $P_{DCS}^{-1}$, $P_{SDTAS_{\tau}}^{-1}$ and $I$($d=d_2\left( x \right)$, $\lambda =1.5$, $\beta =1.2$, $\gamma _1=0.75$)}
		\begin{tabular}{ccccccc}
			\toprule
			\multirow{2}[4]{*}{\textbf{Method}} & \multirow{2}[4]{*}{\textbf{Index}} & \multicolumn{5}{c}{\textbf{N}} \\
			\cmidrule{3-7}          &       & $2^8$   & $2^9 $    & $2^{10}$    & $2^{11}$     &$ 2^{12} $\\
			\midrule
			\multirow{3}[1]{*}{$P_{TAI}^{-1}\left( l \right)$} & $\l$     & 12     & 12    & 12     & 12     & 12 \\
			& IT    & \textbf{12.05 } & \textbf{12.02 } & \textbf{11.01 } & \textbf{11.01 } & \textbf{12.00 } \\
			& CPU   & \textbf{0.35 } & \textbf{0.90 } & \textbf{3.11 } & \textbf{10.89 } & \textbf{60.30 } \\
			\multirow{3}[0]{*}{$P_{CAI}^{-1}\left( l \right)$} & $l$     & 12    & 12    & 12    & 12    & 12 \\
			& IT    & 15.06  & 15.03  & 16.02  & 17.01  & 17.00  \\
			& CPU   & 0.49  & 1.60  & 5.42  & 21.00  & 78.15  \\
			\multirow{3}[0]{*}{$P_{DCS}^{-1}$} & $\beta$     & 0.15  & 0.15  & 0.15  & 0.15  & 0.15 \\
			& IT    & 16.06  & 18.04  & 21.02  & 25.01  & 29.01  \\
			& CPU   & 0.48  & 1.52  & 6.14  & 32.83  & 139.03  \\
			\multirow{3}[0]{*}{$P_{SDTAS_{\tau}}^{-1}$} & $\beta$     & 0.15  & 0.15  & 0.15  & 0.15  & 0.15 \\
			& IT    & 15.06  & 18.04  & 21.02  & 24.01  & 28.01  \\
			& CPU   & 0.24  & 0.65  & 2.84  & 22.00  & 81.56  \\
			\multirow{2}[1]{*}{$I$} & IT    & 125.49 & 174.34 & 233.23 & 300.15 & $\sim$ \\
			& CPU   & 1.46  & 7.41  & 46.45 & 225.74 & $\sim$ \\
			\bottomrule
		\end{tabular}%
		\label{tab:addlabel}%
	\end{table}%
	
	\begin{figure}[!htbp]
		\setlength{\abovecaptionskip}{-0cm}
		\setlength{\belowcaptionskip}{-0.25cm}
		\centering
		\includegraphics[height=5.65cm]{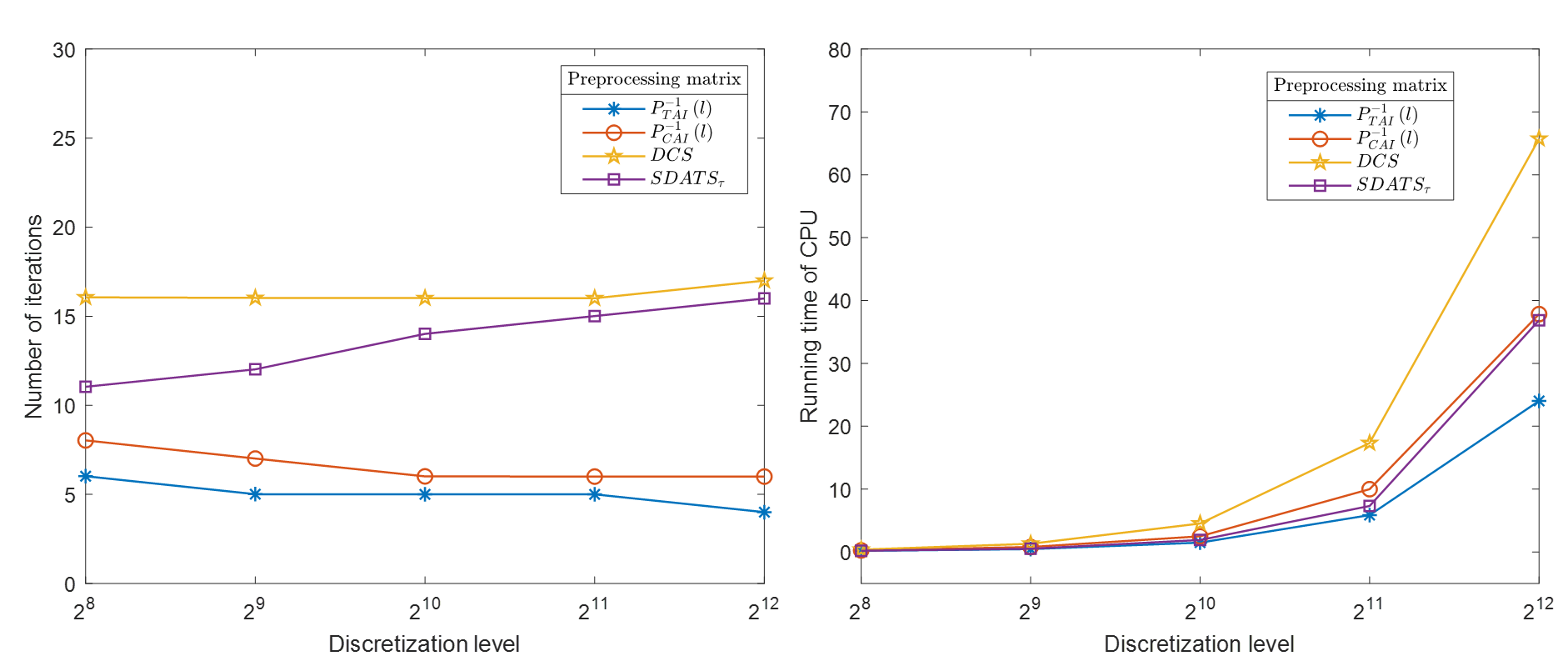}
		\caption*{Fig. 5.5.  Comparison of average iterations and computation time for $P_{TAI}^{-1}\left( l \right)$, $P_{CAI}^{-1}\left( l \right)$, $P_{DCS}^{-1}$, $P_{SDTAS_{\tau}}^{-1}$($d=d_1\left( x \right)$, $\lambda =1.5$, $\beta =1.2$, $\gamma _1=0.75$)}
		\label{figure2}
	\end{figure}
	
	\begin{figure}[!htbp]
		\setlength{\abovecaptionskip}{-0cm}
		\setlength{\belowcaptionskip}{-0.25cm}
		\centering
		\includegraphics[height=7.5cm]{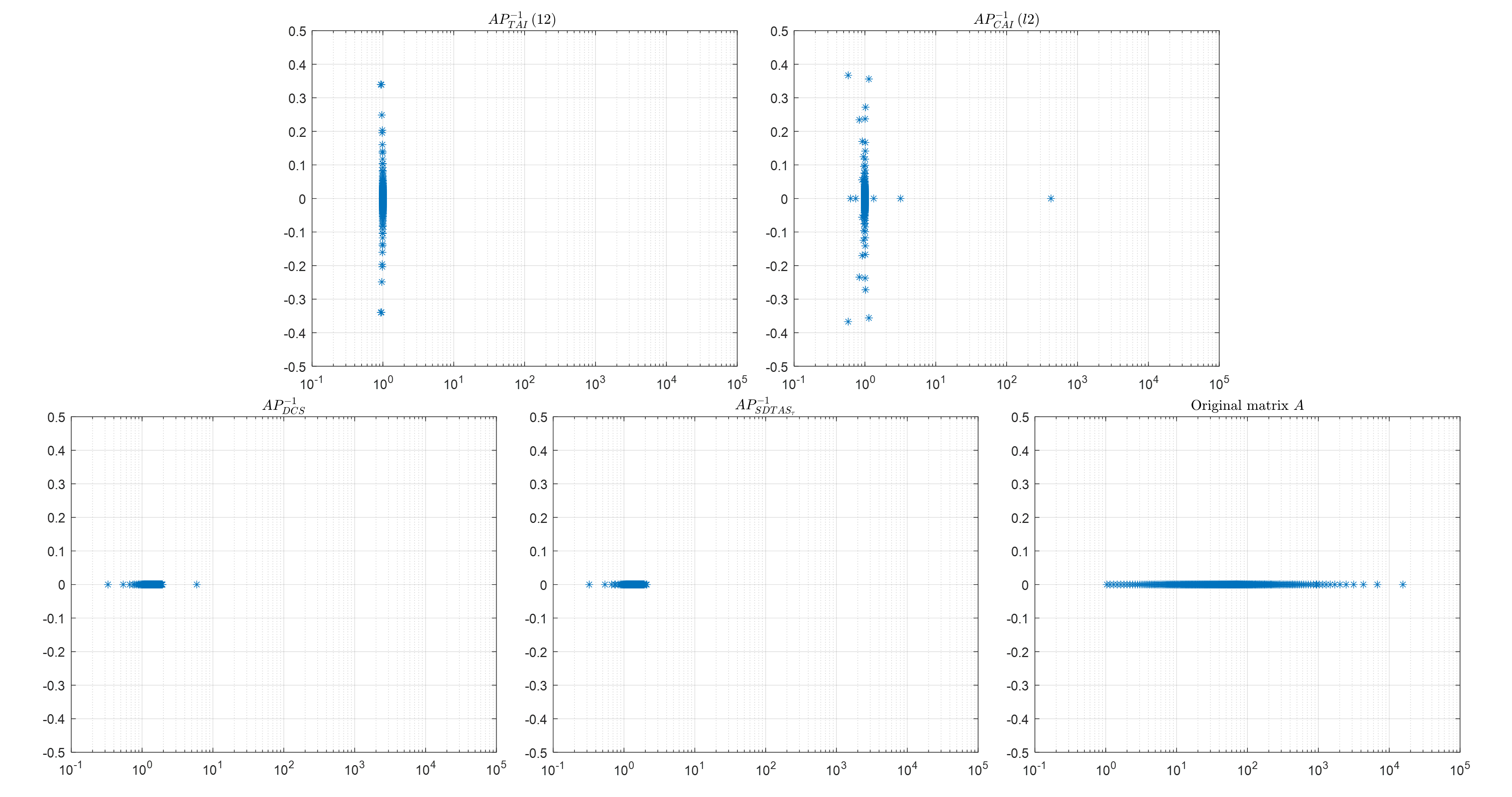}
		\caption*{Fig. 5.6.  The eigenvalue distributions of the preconditioned matrices $AP_{TAI}^{-1}\left( l \right)$, $AP_{CAI}^{-1}\left( l \right)$, $AP_{DCS}^{-1}$, $AP_{SDTAS_{\tau}}^{-1}$ and the original matrix $A$($d=d_2\left( x \right)$, $\lambda =1.5$, $\beta =1.2$, $\gamma _1=0.75$)}
		\label{figure3}
	\end{figure}
	
	We also plot the eigenvalue distributions of the preconditioned coefficient matrices $P_{TAI}^{-1}\left( l \right)$, $P_{CAI}^{-1}\left( l \right)$, $P_{DCS}^{-1}$, $P_{SDTAS_{\tau}}^{-1}$ and $A$ when $N=10^8$, as shown in Fig. 5.6. It can be observed that the eigenvalue distribution of the coefficient matrix without preconditioning is highly dispersed, while the eigenvalue distributions of the preconditioned matrices are concentrated around 1. Notably, the eigenvalue distribution of $P_{TAI}^{-1}\left( l \right)$ is more concentrated than those of the other preconditioners.

	\section{Concluding remarks}\label{section6}
		In this paper, the Crank-Nicolson method and the tempered weighted and shifts Grünwald formula are applied to discretize the tempered fractional diffusion equations. We then get the coefficient matrix of the discretized system having the structure of \( I + DT \), where \( I \) is the identity matrix, \( D \) is a diagonal matrix, and \( T \) is a symmetric positive definite (SPD) Toeplitz matrix. Given the SPD property of the coefficient matrix, we approximate it using \(\tau\) matrix and construct a novel preconditioner using the row-by-row approximate inverse idea proposed by Pan and NG \cite{Pan2014}. The $\tau$ matrix based approximate inverse preconditioner can be efficiently computed using the discrete sine transform(DST). In spectral analysis, the eigenvalues of the preconditioned coefficient matrix are clustered around 1, ensuring fast convergence of Krylov subspace methods with the new preconditioner. Finally, numerical experiments demonstrate the effectiveness of the proposed preconditioner.
	
	\section{Acknowledgement}
		The authors would like to thank the anonymous referees for their constructive comments and suggestions for improving the manuscript.

\end{document}